\documentclass[12pt]{article}  

\usepackage{amsmath}
\usepackage{amssymb}
\usepackage{cite}
\usepackage[usenames]{color}
\usepackage{graphicx}          
\usepackage[dvips]{epsfig}    
\newtheorem{definition}{Definition}

\newtheorem{remark}{Remark}

\newtheorem{lemma}{Lemma}

\newtheorem{proposition}{Proposition}

\title{\LARGE \bf On stabilization of nonlinear systems with drift by time-varying \\ feedback laws
\thanks{
This work was  supported in part by the German Research Foundation (GR 5293/1-1) and
NAS of Ukraine (budget program KPKBK 6541230)\newline
$^{1}$Max Planck Institute for Dynamics of Complex Technical Systems, 39106 Magdeburg, Germany
{\tt\small zuyev@mpi-magdeburg.mpg.de}\newline
$^{2}$Institute of Mathematics, University of W\"{u}rzburg,       97074 W\"{u}rzburg, Germany
        {\tt\small viktoriia.grushkovska@mathematik.uni-wuerzburg.de}
        \newline
$^{3}$Institute of Applied Mathematics and Mechanics, National Academy of Sciences of Ukraine, 841116 Sloviansk, Ukraine
}
}

\author{Alexander Zuyev$^{1,3}$ and Victoria Grushkovskaya$^{2,3}$}
\date{}

\begin{document}

\maketitle
\thispagestyle{empty}

\begin{abstract}
This paper deals with the stabilization problem for nonlinear con\-trol-affine systems with the use of oscillating feedback controls.
We assume that the local controllability around the origin is guaranteed by the rank condition with Lie brackets of length up to~3.
This class of systems includes, in particular, mathematical models of rotating rigid bodies.
We propose an explicit control design scheme with time-varying trigonometric polynomials whose coefficients depend on the state of the system.
The above coefficients are computed in terms of the inversion of the matrix appearing in the controllability condition.
It is shown that the proposed controllers can be used to solve the stabilization problem by exploiting the Chen--Fliess expansion of solutions of the closed-loop system.
We also present results of numerical simulations for controlled Euler's equations and a mathematical model of underwater vehicle to illustrate the efficiency of the obtained controllers.

\end{abstract}

\section{INTRODUCTION}

The stabilization of underactuated mechanical systems with uncontrollable linearization is a challenging mathematical problem related to the theory of critical cases in the sense of Lyapunov.
This problem is of great importance in robotics as the motion of nonholonomic mobile robots and underactuated manipulators is generically described by systems of ordinary differential equations
whose linear approximation does not satisfy Kalman's rank condition at a reference configuration.

On the one hand, it is a well-known fact that even relatively simple models of nonholonomic systems of the form
\begin{equation}
\dot x=\sum_{i=1}^m u_i f_i(x),\quad x\in {\mathbb R}^n,\; u\in {\mathbb R}^m
\label{driftless}
\end{equation}
do not admit ``regular'' stabilizers in the class of time-invariant feedback laws $u=k(x)$ if ${\rm rank}(f_1(0),f_2(0),...,f_m(0))=m<n,$ see~\cite{Bro83}.
On the other hand, each nonlinear control system
\begin{equation}
\dot x = f(x,u),\quad x\in {\mathbb R}^n,\; u\in {\mathbb R}^m,\; f(0,0)=0
\label{general_sys}
\end{equation}
admits a continuous time-varying feedback $u=k(t,x)$ that stabilizes the trivial equilibrium,
provided that small-time local controllability conditions are satisfied at $x=0$ together with some regularity assumptions~\cite{Coron}.
However, the problem of constructing a stabilizing feedback law for an arbitrary controllable system~\eqref{general_sys} (or even~\eqref{driftless} under higher-order controllability condtions)
is far from being satisfactory solved.
It is hardly possible to mention all the achievements in this area due to lack of space,
so we just refer to~\cite{Coron} for a survey of essentially nonlinear design tools.

The goal of this paper is to develop an effective approach for constructing time-varying stabilizers for the class of nonlinear control-affine systems $\dot x =f_0(x)+\sum_{i=1}^m f_i(x)u_i$,
assuming that the local controllability at $x=0$ is ensured by a rank condition with the vector fields $f_1,..., f_m$ and their Lie brackets of the form $[f_i,f_0]$, $[f_i,f_j]$, $[f_i,[f_j,f_{0}]]$, and $[f_i,[f_j,f_k]]$ with some indices $i,j,k=1,2,...,m$.
Our study is motivated by control problems in nonholonomic mechanics and rigid body dynamics, where such type of controllability conditions naturally appears.

The main idea of our construction is summarized in Section~II.
A family of time-periodic feedback controls with coefficients depending on the state vector will be presented explicitly to stabilize the equilibrium of the considered class of systems.
The proposed control algorithm extends our previous design schemes~\cite{ZuSIAM,ZGB16} for the case of systems with drift, i.e. for $f_0(x)\not\equiv 0$.
We will also discuss applications of this control algorithm for the stabilization of a rotating satellite and underwater vehicle in Section~III.

\subsection{Notations, definitions, and auxiliary results}
\begin{definition}
We say that there is a \emph{  resonance of order $N\in\mathbb N$} between pairwise distinct  numbers $k_1,\dots,k_n$, if there exist relatively prime integers  $c_1,{\dots},c_q$  such that $|c_1|+...+|c_q|=N$ and
$c_1k_1+...+c_nk_n=0$.
\end{definition}
\begin{definition}
For a given $\varepsilon>0$, define a partition $\pi_\varepsilon$ of ${\mathbb R}^+=[0,+\infty)$ into the intervals
$
[t_j,t_{j+1})$, $t_j=\varepsilon j$, $j=0,1,2,\dots .$
Given a time-varying feedback law $u=h(t,x)$, $h:{\mathbb R}^+\times D\to\mathbb R^m$, $\varepsilon>0$, and $x^0\in\mathbb R^n$, a \emph{$\pi_\varepsilon$-solution of  system~\eqref{general_sys}}
corresponding to $x^0\in D$ and $h(t,x)$ is an absolutely continuous function  $x(t)\in D$, defined for $t\in {\mathbb R}^+$, such that  $x(0)=x^0$ and
$
\dot x(t)=f\big(x(t), h(t,x(t_j))\big), \quad t\in[t_j,t_{j+1}),
$
for each $j=0,1,2,\,\dots$ .
\end{definition}

 Throughout this paper, $\|a\|$ denotes the Euclidean norm of a vector $a\in\mathbb R^n$, and the norm of a matrix $\cal F$ is defined as $\|{\cal F}\|=\sup\limits_{\|y\|=1}\|{\cal F}y\|$.
For a $\delta>0$, $B_\delta(x^*)$ denotes the $\delta$-neighborhood of $x^*{\in} \mathbb R^n$, and $\overline{B_\delta(x^*)}$ is its closure. For  $h\in C^1(\mathbb R^n;\mathbb R)$ and $\xi\in\mathbb R^n$, 
 $\nabla h(\xi):=\tfrac{\partial h(x)}{\partial x}^\top\Big|_{x=\xi}$. For a function $f:\mathbb R\to\mathbb R$,  $f(z)=O(z)$ as $z\to 0$ means that there is a $c>0$ such that $|f(z)|\le c|z|$ in some neighborhood of $0$.
   For  $f,g:\mathbb R^n\to\mathbb R^n $, the directional derivative is denoted as
 $ L_gf(x)=\tfrac{\partial f(x)}{\partial x}g(x)$, and  $[f,g](x)= L_fg(x)- L_gf(x)$ stands for the Lie bracket.
 For any set $S$ and $n\in\mathbb N$, $S^n$ denotes the $n$-th Cartesian power of $S$, i.e.  $S^n=\underbrace{S\times S\times\dots\times S}_{n}=\{(x_1,\dots,x_n):x_i\in S \text{ for every } i\in\{1,\dots,n\}\}$.


\begin{lemma}[\cite{La95,GZE18}]\label{volterra}
{ Let the vector fields $f_i$ be Lipschitz continuous in a domain $D\subseteq\mathbb R^n$, and $f_i\in C^{\nu+1}(D\setminus\Xi;\mathbb R)$, where $\Xi=\{ x\in D:f_i( x)=0\text{ for all }1\le i\le \ell\}$, $\nu\ge 1$. Assume, moreover, that $ L_{f_j}f_i, L_{f_l}L_{f_j}f_i\in
C(D;\mathbb R^n)$,   for all $i,j,l=\overline{1,\ell}$. Let $ x(t)\in D$, $t\in[0,\tau]$, be a  solution of the system
$\displaystyle\dot  x=\sum_{i=1}^\ell f_i( x)w_i(t)$
 with $u\in C([0,\tau];\mathbb R^m)$
 and $x(0)=x^0\in D$. Then, for any $t\in[0,\tau]$,  $ x(t)$ can be represented by the $\nu$-th order Chen--Fliess series:}
\begin{align}
     x(t)&=x^0+\sum_{i=1}^{ \nu}\sum_{j_1,\dots,j_ i=1}^\ell L_{f_{j_{ i}}}L_{f_{j_{ i-1}}}\dots f_{j_1}( x^0)\label{volt1}\\
     &\times\int_0^t\int_0^{s_{1}}\dots\int_0^{s_{ i-1}} w_{j_1}(s_1)\dots w_{j_{ i}}(s_ i)ds_{ i}\dots ds_1 +r(t),\nonumber
   \end{align}
 where the remainder is
$$
\begin{aligned}
r(t)=\sum_{j_1,\dots,j_{\nu+1}=1}^\ell&\int_0^t\int_0^{s_{1}}\dots\int_0^{s_{\nu}} L_{f_{j_{\nu+1}}}L_{f_{j_{\nu}}}\dots f_{j_1}(x(s_{\nu+1}))\\
     &\times w_{i_1}(s_1)\dots w_{j_{\nu+1}}(s_{\nu+1})ds_{\nu+1}\dots ds_1.
\end{aligned}
$$
     \end{lemma}
     Note that the above Chen--Fliess expansion can be obtained from the Volterra series~\cite{La95}.
\section{Main results}
In this section, we consider a class of nonlinear control systems of the form
\begin{equation}\label{nonh}
\begin{aligned}
&\dot x =f_0(x)+\sum_{k=1}^mf_k(x)u_k\\
\end{aligned}
\end{equation}
where  $x=(x_1,\dots,x_n)^\top\in D$ is the state vector, $D\subset \mathbb R^n$ is a domain, $0\in D$,
 $u=(u_1,\dots,u_m)^\top$ is the control, and the vector fields $f_0,f_1,\dots,f_m\in C^2(\mathbb R^n;\mathbb R^n)$ satisfy the following rank condition: for all $x\in D$,
 \begin{equation}\label{rank}
\begin{aligned}
{\rm span}\Big\{ f_i(x),\, [f_{i_1},f_{i_2}]&(x),\,\big[f_{j_1},[f_{j_2},f_{j_3}]\big](x),\\
& [f_l,f_0](x),\,\big[f_{l_1},[f_{l_2},f_{0}]\big](x) \Big\}={\mathbb R}^n,
\end{aligned}
\end{equation}
where $i\in S_1$, $(i_1,i_2)\in S_2$, $(j_1,j_2,j_3)\in S_3$, $l\in S_{10}$, $(l_1,l_2)\in S_{20}$,
with some set of indices $S_1,S_{10}\subseteq \{1,2,\dots,m\}$, $S_2,S_{20}\subseteq \{1,2,\dots,m\}^2$, $S_3\subseteq \{1,2,\dots,m\}^3$ such that $|S_1|+|S_2|+|S_3|+|S_{10}|+|S_{20}|=n$.
Condition~\eqref{rank} can be used for checking the local controllability of system~\eqref{nonh} at the origin~\cite{Sussmann87,AgrachevSarychev}.
In~\cite{ZG17}, we have considered local approximate steering and local path-following problems for the case $S_2=\emptyset$, $S_3=\emptyset$, and $S_{20}=\emptyset$. However, the results of the above paper do not guarantee any stability properties. Below we propose  novel stabilizability results for system~\eqref{nonh} satisfying~\eqref{rank} and introduce explicit formulas for the control design.
To simplify the presentation, we  first consider two cases:

1) $|S_1|+|S_{10}|=n$, i.e. the vector fields of the system satisfy
\begin{equation}\label{rank1}
\begin{aligned}
{\rm span}\left\{ f_i(x),\, [f_l,f_0](x) \right\}={\mathbb R}^n;
\end{aligned}
\end{equation}

2) $|S_1|+|S_{20}|=n$, i.e. the vector fields of system satisfy
\begin{equation}\label{rank2}
\begin{aligned}
{\rm span}\Big\{f_i(x),\,\big[f_{l_1},[f_{l_2},f_{0}]\big](x) \Big\}={\mathbb R}^n.
\end{aligned}
\end{equation}

Then we will show how to combine the above particular cases with the controls from~\cite{GZ18} in case when the rank condition~\eqref{rank} is satisfied.
Note that the rank condition~\eqref{rank2} and its iterations appears as controllability assumptions in solid and fluid mechanics for affine systems controlled by forces and torques~\cite{AgrachevSarychev,BulloLewis}.

To stabilize system~\eqref{nonh} at $x=0$, we use a parameterized family of control functions $u_k=u_k^\varepsilon(t,x)$ consisting of several terms, each of which is aimed to implement the motion in the direction of a certain vector field from the rank condition for small values of the parameter $\varepsilon>0$.
\subsection{Stabilization of system~\eqref{nonh} under condition~\eqref{rank1}}
For the case 1), our proposed control design is as follows:
\begin{equation}\label{cont1}
\begin{aligned}
u_k^\varepsilon(t,x) = &\sum_{i\in S_1}h^k_i(x)+\tfrac{1}{\varepsilon}\sum_{l\in S_{10}} h^k_{l0}(t,x)\\
\end{aligned}
\end{equation}
with $h^k_i(x) = \delta_{ki}a_i(x)$,
\begin{equation*}
h^k_{l0}(t,x)=\delta_{kl}{2\pi \kappa_{l0}a_{l0}(x)}\sin \Big(\tfrac{2\pi \kappa_{l0} t}{\varepsilon}\Big),
\quad k=1,\dots,m,
\end{equation*}
where $\delta_{ij}$ is the Kronecker delta,  $\kappa_{l0}$ are pairwise distinct  positive integers. Denote the vector   of real-valued coefficient functions
$a_i(x)$ and $a_{l0}(x)$ as $a^{(1)}(x)=\Big((a_{i})_{i\in S_1}\ \ (a_{l0})_{l\in S_{10}}\Big)^\top\in\mathbb R^n$. The goal is to define it in such a way that  a certain potential function $V\in C^2(\mathbb R^n;\mathbb R^+)$ decreases along the trajectories of system~\eqref{nonh}. With this purpose, we put
\begin{equation}\label{a1}
\begin{aligned}
a^{(1)}(x)=-\mathcal F_1^{-1}(x)\big(\gamma\nabla V(x)+f_0(x)\big),
\end{aligned}
\end{equation}
where $\gamma>0$ plays the role of control gain, and $\mathcal F_1^{-1}(x)$ is the inverse matrix for
$  \mathcal F_1(x)= \Big(\big(f_{i}(x)\big)_{i\in S_1}\  \big([f_{l},f_{0}](x)\big)_{l\in S_{10}}\Big).
$
Obviously, $\mathcal F_1^{-1}(x)$ exists whenever the rank condition~\eqref{rank} holds.
Note that the proposed control formulas are much simpler compared to the ones used in~\cite{ZG17}.

The main idea behind our control design approach is the approximation of  trajectories of the auxiliary system
$
  \dot{\bar x}=-\gamma \nabla V(\bar x)$, $ \bar x\in D,\, \bar x(0)=x(0)
$
by the trajectories of~\eqref{nonh}.
Indeed, under the proposed choice of control functions, the representation~\eqref{volt1} yields
\begin{equation}\label{volt_eps_1}
x(\varepsilon)=x^0+\varepsilon\big(f_0(x^0)+\mathcal F_1(x^0)a^{(1)}(x^0)\big)
+\Omega_1(x^0,\varepsilon)+R(\varepsilon),
\end{equation}
where  $\Omega_1(x^0,\varepsilon)$ is given in the Appendix, and the remainder $R_1(\varepsilon)$ is calculated as $r(\varepsilon)$ in Lemma~\ref{volterra} with the summation indices $j_1,j_2,j_3\in\{0,1,\dots,m\}$ and $w_0(t)\equiv 1$. Substituting~\eqref{a1} into~\eqref{volt_eps_1}, we obtain
$
x(\varepsilon)=x^0-\gamma\varepsilon\nabla V(x^0)+\Omega_1(x^0,\varepsilon)+R(\varepsilon).
$
Let us take $V(x)=\tfrac{1}{2}\|x\|^2$ 
and suppose that
$$
f_0(0)=L_{f_0}f_0(0)=0\eqno{(A1)}.
$$
Since $f_0\in C^2(\mathbb R^n;\mathbb R^n)$, (A1) means that $\|f_0(x)\|$ and $\|L_{f_0}f_0(x)\|$ are $O(\|x\|)$ as $x\to0$.
Thus, for each compact set $D_0\subseteq D$, $0{\in} D_0$, there exist  $\sigma_{i}\ge 0$ such that, for any $x^0{\in} D_0$,
$$
\begin{aligned}
\|x(\varepsilon)\|\le \|x^0\|\big(1-\varepsilon\gamma\big)+&\varepsilon\sigma_1\|x^0\|^2+\varepsilon^2\sigma_2\|x^0\|+\|R_1(\varepsilon)\|.
\end{aligned}
$$
Based on this estimate, the following result can be proved.
\begin{proposition}\label{thm_s10}
  \emph{Let $D\subseteq \mathbb R^n$, $f_i\in C^3(D;\mathbb R^n)$, $i=0,\dots,m$. Suppose that assumption~(A1) holds and, furthermore,   there exists an  $\alpha>0$ such that
 $ \|\mathcal F_1^{-1}(x)\|\le \alpha \text{ for all }x\in D$.
If the functions $u_k=u_k^\varepsilon (t,x)$, $k=1,\dots,m$, are defined by~\eqref{cont1}--\eqref{a1}, then
 there exist  $\gamma,\delta,\bar\varepsilon>0$  such that, for any $\varepsilon\in(0,\bar\varepsilon]$, each $\pi_\varepsilon$-solution of system~\eqref{nonh} with the initial data $x(0)=x^0\in B_\delta(0)$ is well-defined on $t\in {\mathbb R}^+$ and
$
\|x(t)\|\to 0 \text{ as }t\to\infty.
$
}
\end{proposition}

The proof of the above proposition is based on subtle estimates of the remainder $R_1(\varepsilon)$ and the techniques from~\cite{ZuSIAM,GZ18}. In particular, it goes along the same line as the proofs of~\cite[Theorem~2.2]{ZuSIAM} and~\cite[Theorem~1]{GZ18} with taking into account the  drift term $f_0(x)$.
The detailed proof of Proposition~\ref{thm_s10}, including description of connections between $\varepsilon,\gamma,\delta$, will be presented in the extended version of the paper.
\begin{remark}~\label{rem_non0}
  If the vector fields of~\eqref{nonh} fail to satisfy~(A1), a weaker result on practical stabilizability of system~\eqref{nonh} can be deduced using the techniques proposed, e.g., in~\cite{GDEZ17}. We will illustrate this case with an example in Section~III.B.
\end{remark}
\subsection{Stabilization of system~\eqref{nonh} under condition~\eqref{rank2}}
For the case 2), we propose the following controls:
\begin{equation}\label{cont2}
\begin{aligned}
u_k^\varepsilon(t,x) = &\sum_{i\in S_1}h^k_i(x)+\tfrac{1}{\varepsilon}\sum_{(l_1,l_2)\in S_{20}} h^k_{l_1l_20}(t,x),
\end{aligned}
\end{equation}
where $h^k_i(x) =\delta_{ki} a_i(x)$,
\begin{equation*}
\begin{aligned}
h^k_{l_1l_20}(t,x)= &4\pi\kappa_{l_1l_20}\sqrt{|a_{l_1l_20}(x)|}\cos\Big(\tfrac{2\pi \kappa_{l_1l_20}t}{\varepsilon}\Big)\\
&\quad\times\big(\delta_{kl_1}+\delta_{kl_2}{\rm sign}\big(a_{l_1l_20}(x)\big)\big),
\end{aligned}
\end{equation*}
for $k=1,\dots,m$, and  pairwise distinct  positive integers $\kappa_{l_1l_20}$. Moreover, we assume that there are no resonances of order 2 between $\kappa_{l_1l_20}$. Similarly to the previous subsection, we define the  vector $a^{(2)}(x)=\Big((a_{i})_{i\in S_1}\ (a_{l_1l_20})_{(l_1,l_2)\in S_{20}}\Big)^\top$
as
\begin{equation}\label{a2}
a^{(2)}(x)=-\mathcal F_2^{-1}(x)\big(\gamma\nabla V(x)+f_0(x)\big),\quad\gamma>0,
\end{equation}
where  $\mathcal F_2^{-1}(x)$ is the inverse for the matrix
$$
\Big(\big(f_{i}(x)\big)_{i\in S_1}\  \big(\big[f_{l_1},[f_{l_2},f_0]\big](x){+}\big[f_{l_2},[f_{l_1},f_0]\big](x)\big)_{l\in S_{20}}\Big).
$$
In this case, the expansion~\eqref{volt1} takes the form
\begin{equation}\label{volt_eps_2}
x(\varepsilon)=x^0-\gamma\varepsilon\nabla V(x^0)+\Omega_2(x^0,\varepsilon)+R_2(\varepsilon).
\end{equation}
The explicit formula for $\Omega_2(x^0,\varepsilon)$ is given in the Appendix, and the remainder $R_2(\varepsilon)$ is calculated according to Lemma~\ref{volterra}.
To ensure that the function $V(x)=\tfrac{1}{2}\|x\|^2$ decays over the time period $\varepsilon$ along the trajectories of system~\eqref{nonh} with controls~\eqref{cont2}--\eqref{a2}, we suppose that assumption (A1) holds together with the following properties:
$$
\begin{aligned}
&L_{f_0}L_{f_0}f_0(0)=\big[f_{0},[f_{0},f_{k}]\big](0)=0,\\
&\big[f_{l_1},[f_{l_1},f_0]\big](x)+\big[f_{l_2},[f_{l_2},f_0]\big](x)=O(\|x\|^\mu),\qquad\qquad(A2)\\
&\big[f_{0},[f_{l_1},f_{l_2}]\big](x)=O(\|x\|^\mu)\;\text{as }\|x\|\to 0\;\text{ with some }\mu>0,\\
&\text{for any}\,(l_1,l_2)\in S_{20}\, \text{ and any}\;k:  (l_1,k)\in S_{20} \text{ or }(k,l_2)\in S_{20}.
\end{aligned}
$$
Then, for each compact set $D_0\subseteq D$, $0\in D_0$, there exist constants $\tilde\sigma_{i}\ge 0$ such that, for any $x^0\in D_0$,
$$
\|x(\varepsilon)\|{\le} \|x^0\|\big(1-\varepsilon\gamma\big)+\varepsilon\tilde\sigma_1\|x^0\|^{1+\mu}+\varepsilon^2\tilde\sigma_2\|x^0\|+\|R_2(\varepsilon)\|.
$$
Again, estimating the remainder $R_2(\varepsilon)$ and using the techniques from~\cite{ZuSIAM,GZ18}, we can state the following result.

\begin{proposition}\label{thm_s20}
  \emph{Let $D\subseteq \mathbb R^n$, $f_i\in C^4(D;\mathbb R^n)$, $i=0,\dots,m$. Suppose that assumptions~(A1) and~(A2) hold and, furthermore,   there exists an  $\alpha>0$ such that
 $ \|{\mathcal F_2}^{-1}(x)\|\le \alpha \text{ for all }x\in D$.
If the functions $u_k=u_k^\varepsilon (t,x)$, $k=1,\dots,m$, are defined by~\eqref{cont2}--\eqref{a2}, then
 there exist  $\gamma,\delta,\bar\varepsilon>0$  such that, for any $\varepsilon\in(0,\bar\varepsilon]$, each $\pi_\varepsilon$-solution of system~\eqref{nonh} with the initial data $x(0)=x^0\in B_\delta(0)$ is well-defined on $t\in {\mathbb R}^+$ and
$
\|x(t)\|\to 0 \text{ as }t\to\infty.
$
}
\end{proposition}
\subsection{Stabilization of system~\eqref{nonh} under condition~\eqref{rank}}
If  the rank condition~\eqref{rank} involves Lie brackets of the type $ [f_{i_1},f_{i_2}]$, $\big[f_{j_1},[f_{j_2},f_{j_3}\big]$,
$[f_l,f_0]$, and $\big[f_{l_1},[f_{l_2},f_{0}]\big]$, then stabilizing controllers can be constructed on the basis of formulas~\eqref{cont1}, \eqref{cont2}, and the scheme from~\cite{GZ18} as follows:
\begin{equation}\label{cont}
\begin{aligned}
u_k^\varepsilon(t,x) &= \sum_{i\in S_1}h^k_i(x)+\tfrac{1}{\varepsilon}\sum_{l\in S_{10}} h^k_{l0}(t,x)+\tfrac{1}{\sqrt\varepsilon}\sum_{(i_1,i_2)\in S_{2}} h^k_{i_1i_2}(t,x)\\
&+\tfrac{1}{\varepsilon}\sum_{(l_1,l_2)\in S_{20}} h^k_{l_1l_20}(t,x)+\tfrac{1}{\sqrt[3]{\varepsilon^2}}\sum_{(j_1,j_2,j_3)\in S_{3}} h^k_{j_1j_2j_3}(t,x),
\end{aligned}
\end{equation}
where $h^k_i(x)$, $ h^k_{l0}(t,x)$, $h^k_{l_1l_20}(t,x)$ are constructed as in Sections~II.A, II.B, and
\begin{equation*}
\begin{aligned}
h^k_{i_1i_2}&(t,x) =2\sqrt{\pi \kappa_{i_1i_2}|a_{i_1i_2}|}\Big(\delta_{ki_1}{\rm sign}(a_{i_1i_2}(x))\cos\Big(\tfrac{2\pi \kappa_{i_1i_2}t}{\varepsilon}\Big)\\
 &+\delta_{ki_2}\sin\Big(\tfrac{2\pi \kappa_{i_1i_2}t}{\varepsilon}\Big)\Big), \\
h^k_{j_1j_2j_3}&(t,x) =2\sqrt[3]{2\pi^2 (\kappa_{1 j_1 j_2 j_3}^2-\kappa_{2 j_1 j_2 j_3}^2)a_{j_1j_2j_3}(x)}\\
&\times\bigg(\sin\Big(\tfrac{2\pi \kappa_{1 j_1 j_2 j_3}t}{\varepsilon}\Big)\Big(\delta_{kj_1}\cos\Big(\tfrac{2\pi \kappa_{2 j_1 j_2 j_3}t}{\varepsilon}\Big)+\delta_{k j_2}\Big)\\
&+ \delta_{k j_3}\cos\Big(\tfrac{2\pi \kappa_{2 j_1 j_2 j_3}t}{\varepsilon}\Big)\bigg).
\end{aligned}
\end{equation*}
In this case, the conditions of Propositions~\ref{thm_s10} and~\ref{thm_s20} should be satisfied with the corresponding components $h^k_i(x)$, $ h^k_{l0}(t,x)$, and $h^k_{l_1l_20}(t,x)$. The integers $\kappa_{l0}$, $\kappa_{l_1l_20}$, $\kappa_{i_1i_2}$, $\kappa_{1 j_1 j_2 j_3}$, $\kappa_{2 j_1 j_2 j_3}$, $\kappa_{3 j_1 j_2 j_3}=\kappa_{1 j_1 j_2 j_3}+\kappa_{2 j_1 j_2 j_3}$, $\kappa_{4 j_1 j_2 j_3}=\kappa_{1 j_1 j_2 j_3}-\kappa_{2 j_1 j_2 j_3}$ have to be
positive and pairwise distinct. Moreover, if $S_{20}\ne \emptyset$ or $S_{3}\ne\emptyset$, then we assume
 that there are no resonances of order 2 between the above-listed frequencies (see Def.~1). This assumption is needed to reduce the number of Lie brackets excited by controls~\eqref{cont}.
  The  vector of state-dependent coefficients
 $$
 \begin{aligned}
 a(x)=\Big((a_{i})_{i\in S_1}\ (a_{i_1i_2})&_{(i_1,i_2)\in S_{2}}\ (a_{j_1j_2j_3})_{(j_1,j_2,j_3)\in S_{3}}\\
 & (a_{l0})_{l\in S_{10}}\ (a_{l_1l_20})_{(l_1,l_2)\in S_{20}}\Big)^\top
\end{aligned}
$$
is defined as
$
a(x)=-\mathcal F^{-1}(x)\big(\gamma\nabla V(x)+f_0(x)\big),
$
where $\gamma>0$, and $\mathcal F^{-1}(x)$ is the inverse matrix for
$$
\begin{aligned}
&\Big(\big(f_{i}(x)\big)_{i\in S_1}\ ([f_{i_1},f_{i_2}])_{(i_1,i_2)\in S_{2}}\ \big(\big[f_{j_1},[f_{j_2},f_{j_3}]\big]\big)(x)_{(j_1,j_2,j_3){\in} S_{3}}\\
&\ \big([f_{l},f_{0}](x)\big)_{l\in S_{10}}\ \big(\big[f_{l_1},[f_{l_2},f_0]\big](x){+}\big[f_{l_2},[f_{l_1},f_0]\big](x)\big)_{l\in S_{20}}\Big).
\end{aligned}
$$
Stabilizability conditions for this case will be formulated in the extended version of this work.
\section{EXAMPLES}
\begin{figure*}[tpt]
 \begin{minipage}{0.48\linewidth}
\begin{center}
\includegraphics[width=1\linewidth]{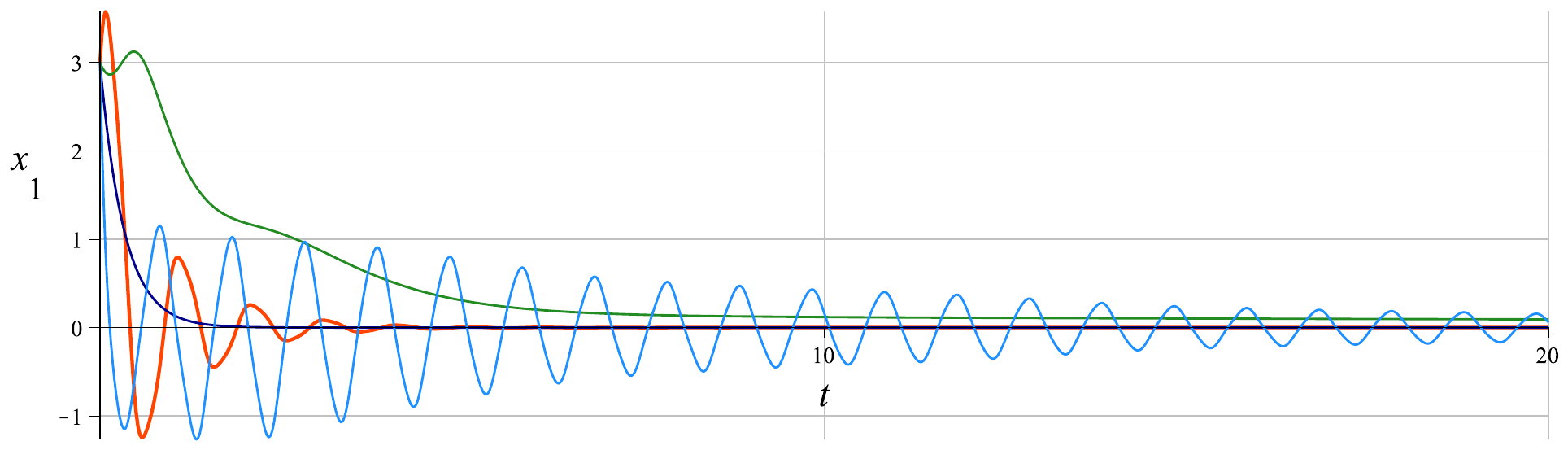}
\includegraphics[width=1\linewidth]{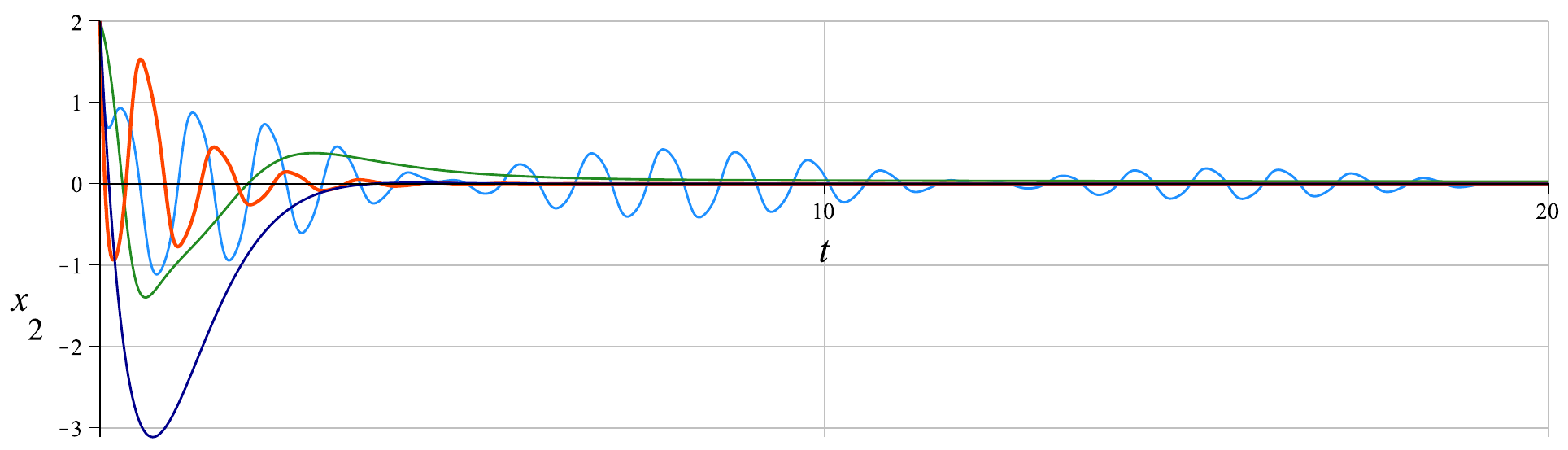}
\includegraphics[width=1\linewidth]{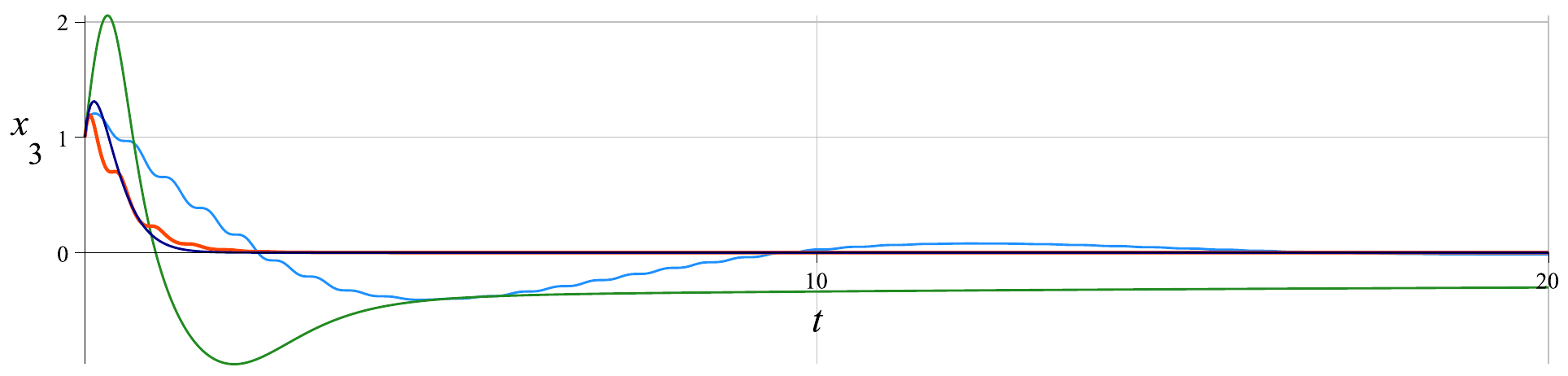}
a)
\end{center}
 \end{minipage}\hfill
 \begin{minipage}{0.48\linewidth}
\begin{center}
\includegraphics[width=1\linewidth]{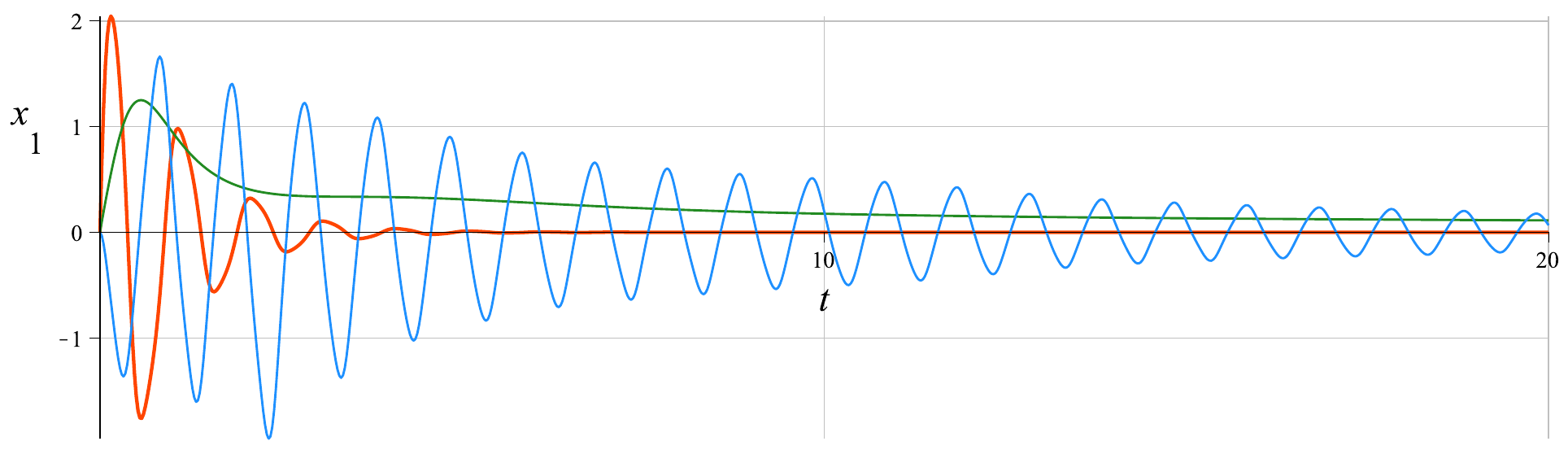}
\includegraphics[width=1\linewidth]{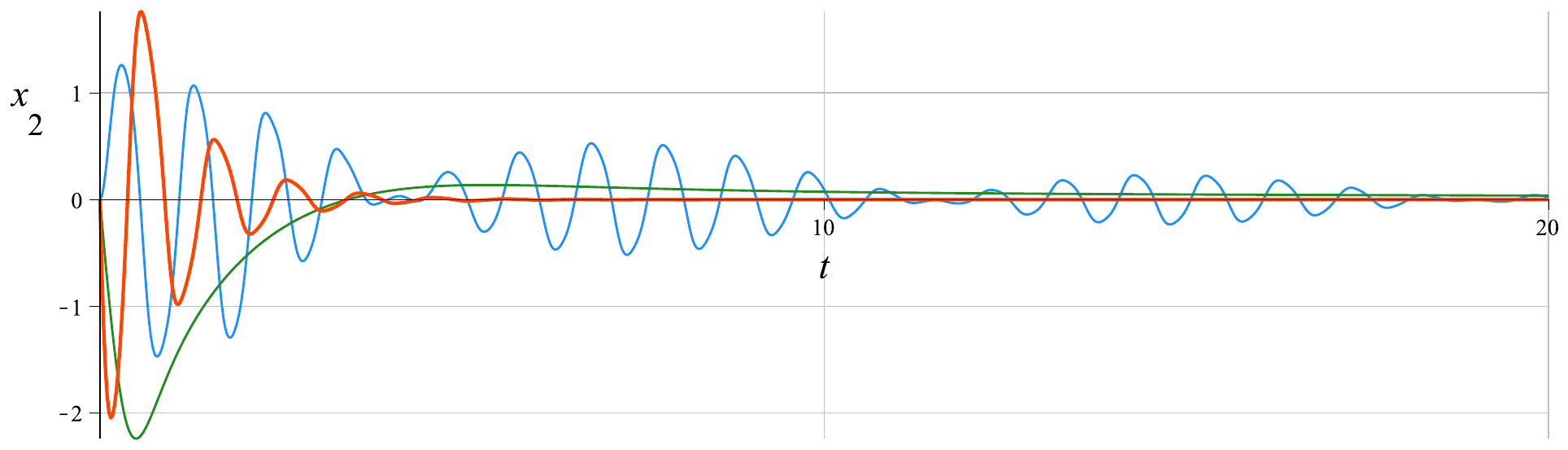}
\includegraphics[width=1\linewidth]{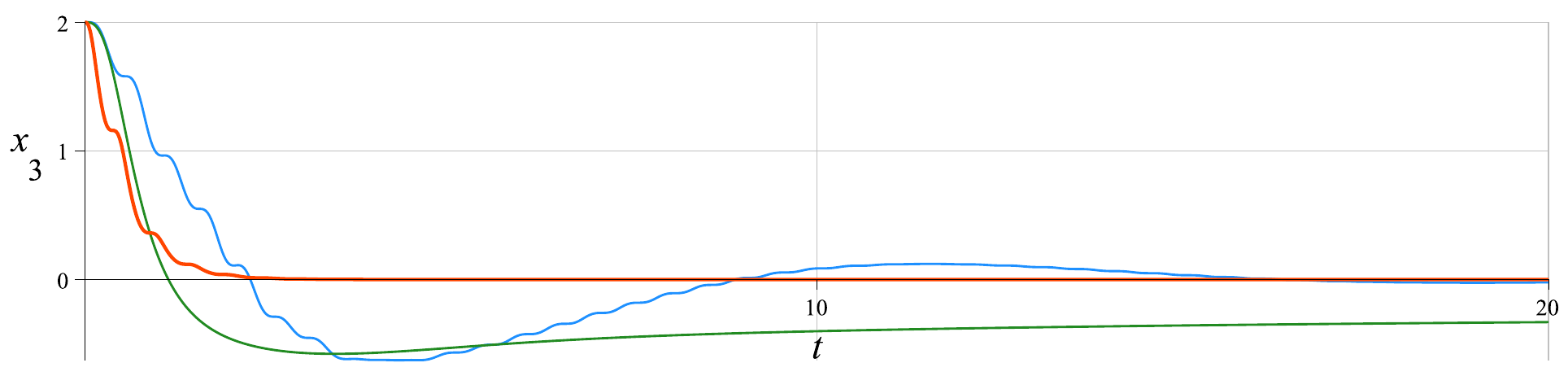}
b)
\end{center}
 \end{minipage}
\begin{center} \begin{minipage}{0.75\linewidth}
\includegraphics[width=1\linewidth]{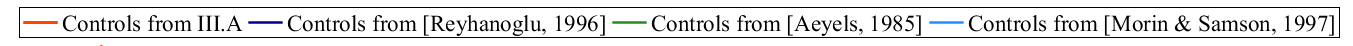}
\end{minipage}
\end{center}
\caption{Time-plots of the trajectories of system~\eqref{Euler} with controls~\eqref{Euler-feedback} (red) and controls from~\cite{Aeyels} (green), \cite{Rey96} (dark blue), and ~\cite{Mor97} (light blue). Figure~a):  $x(0)=(3,2,1)^\top$; figure~b): $x(0)=(0,0,2)^\top$.
\label{fig_body}}
\hspace{-1.5em}
\begin{minipage}{0.5\linewidth}
\begin{center}
\includegraphics[width=0.8\linewidth]{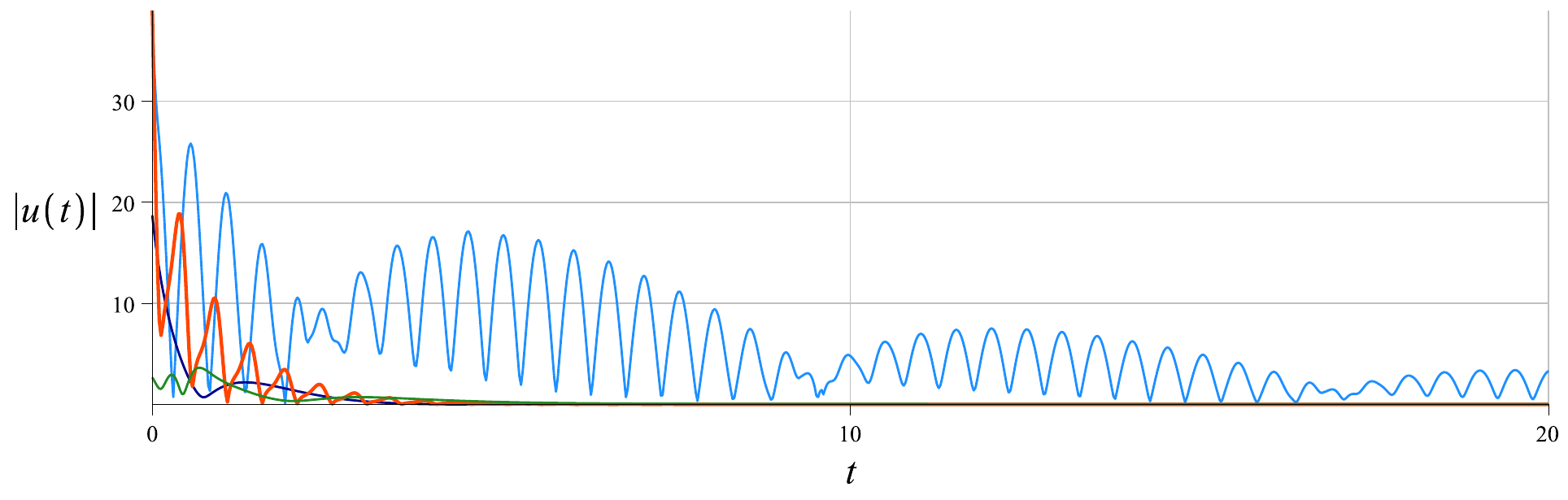}
\end{center}
\end{minipage}
\begin{minipage}{0.5\linewidth}
\begin{center}
\includegraphics[width=0.8\linewidth]{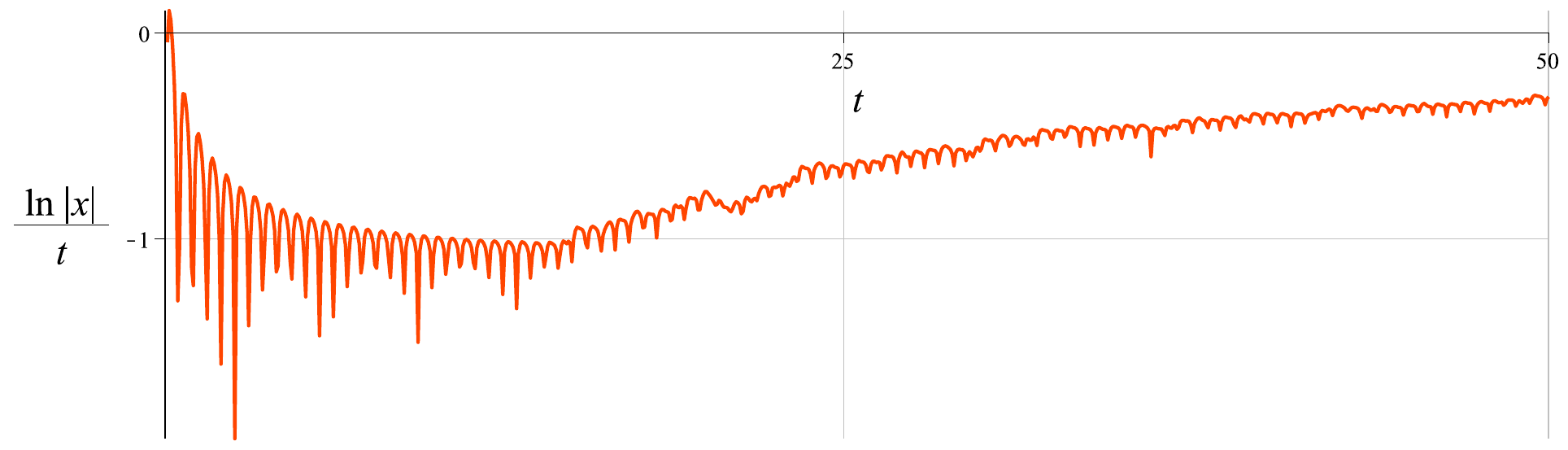}
\end{center}
\end{minipage}
\caption{Left: time-plot of the function $\|u(t)\|$, where $u(t)$ is the control given by~\eqref{Euler-feedback} (red) and from~\cite{Aeyels} (green), \cite{Rey96} (dark blue),~\cite{Mor97} (light blue). Right: time-plot of the function $\tfrac{1}{t}\ln\|x(t)\|$, where $x(t)$ is the solution of system~\eqref{Euler} with controls~\eqref{Euler-feedback}. In both figures, $x(0)=(3,2,1)^\top$.
\label{fig_cont}}
\end{figure*}
\subsection{Stabilization of a rotating rigid body}
Consider Euler's equations for a rigid body rotating around its center of mass:
\begin{equation}
 \dot x_1 = \alpha_1 x_2 x_3 + u_1,\; \dot x_2 = \alpha_2 x_1 x_3 + u_2,\; \dot x_3 = \alpha_3 x_1 x_2,\\
\label{Euler}
\end{equation}
where $x_1,x_2,x_3$ are projections of the angular velocity vector on the principal axes of inertia of the body,
and $u_1,u_2$ correspond to the control torques with respect to the first and the second principal axis, respectively.
The parameters $\alpha_i$ are related to the central moments of inertia of the rigid body $J_1,J_2,J_3$ as
$
\alpha_1 = \tfrac{J_2-J_3}{J_1},\;   \alpha_2 = \tfrac{J_3-J_1}{J_2},\; \alpha_3 = \tfrac{J_1-J_2}{J_3}.
$
In the sequel, we assume that $\alpha_3\neq 0$.
System~\eqref{Euler} represents the angular motion of a spacecraft as an absolutely rigid body without moving masses.
The control torques $u_1$ and $u_2$ can be generated by jet engines, and the change of mass due to the operation of engines is neglected.

The stabilization problem for Euler's equation with time-invariant feedback controls has been already thoroughly studied by several authors.
The orientation of a satellite along a given direction (uniaxial stabilization) was considered in~\cite{Zubov}
for Euler's equations coupled with Poisson's equations using a three-dimensional control (fully actuated case).
It was also shown in~\cite{Zubov} that the stabilization of two given directions is possible in the fully actuated case.
A single-axis stabilization problem for the underactuated Euler--Poisson equations was solved in~\cite{Zu2001}.

By applying the feedback transformation
$y_1=x_1$, $y_2=x_2$, $y_3=x_3/\alpha_3$, $w_1 = \alpha_1 x_2 x_3 + u_1$, $w_2=\alpha_2 x_1 x_3 + u_2$, system~\eqref{Euler} takes the form
\begin{equation}
\dot y_1 = w_1,\; \dot y_2 = w_2,\; \dot y_3 = y_1 y_2.
\label{Euler-transformed}
\end{equation}
As it was proved in~\cite{Bro83} and~\cite{Aeyels}, the equilibrium $y=0$ of~\eqref{Euler-transformed} is asymptotically stabilizable by a smooth time-invariant feedback law.
However, the analysis in~\cite{Aeyels} is based on the center manifold approach, and the resulting reduced system $\dot z = -\varkappa z^5 + o(|z|^5)$ does not exhibit exponential decay rate as $t\to +\infty$.
The trivial equilibrium of system~\eqref{Euler-transformed} is shown to be stabilizable in finite time by means of discontinuous state feedback controls~\cite{Jammazi}.
The problem of robust stabilization of system~\eqref{Euler} with external disturbances was addressed in~\cite{Astolfi-Rapaport}.

Note that the trivial equilibrium of Euler's equations is stabilizable by a one-dimensional control acting along a ``skewed'' direction, provided that the rigid body is asymmetric~\cite{Aeyels-Szafranski}.
This result was extended in~\cite{Sontag-Sussmann} to the body with two identical principal moments of inertia.
It was also shown there that the body with a spherical tensor of inertia cannot be stabilized by a one-dimensional control.

In this section, we will illustrate the behavior of trajectories of system~\eqref{Euler} with controls~\eqref{cont2}.
We will present simulation results to show that the solutions of the closed-loop system decay exponentially with time.

Let us rewrite control system~\eqref{Euler} in the vector form~\eqref{nonh}:
\begin{equation}
\dot x = f_0(x) + u_1 f_1(x) + u_2 f_2(x),\quad x\in{\mathbb R}^3,\;u\in{\mathbb R}^2,
\label{Euler-vector}
\end{equation}
with
$
f_0 = \big(
\alpha_1 x_2 x_3,
\alpha_2 x_1 x_3,
\alpha_3 x_1 x_2,
\big)^\top,\;
f_1 = \big(
1, 0,\ 0\big)^\top,\;
f_2 = \big(
0, 1, 0\big)^\top.
$
System~\eqref{Euler-vector} satisfies the controllability condition of the type~\eqref{rank2}:
$$
{\rm span}\{f_1,f_2,[f_1,[f_2,f_0]] \} = {\mathbb R}^3\;\text{for all}\; x\in{\mathbb R}^3,
$$
provided that $\alpha_3\neq 0$, i.e. $J_1 \neq J_2$. Thus, it is easy to see that the assumptions (A1) and (A2) are satisfied. In particular,
$\big[f_1,[f_1,f_0]\big](x){\equiv} 0$, $\big[f_2,[f_2,f_0]\big](x){\equiv} 0$, $\big[f_0,[f_1,f_2]\big](x){\equiv }0$,
and
$$\big[f_0,[f_0,f_j]\big](x)=O(\|x\|^2)\;\;\text{for} \; \|x\|\to 0, \quad j=1,2.$$
Furthermore,
the matrix
$
{\mathcal F_2}^{-1}(x)= \begin{pmatrix}1 & 0 & 0 \\ 0 & 1 & 0 \\ 0 & 0 & \tfrac{1}{2\alpha_3} \end{pmatrix}
$
is obviously nonsingular and satisfies the conditions of Proposition~\ref{thm_s20} with $\alpha=2+\tfrac{1}{2\alpha_2}$.
Thus, we can apply the stabilization scheme proposed in Section~2.B with $S_1=\{1,2\}$ and $S_{20}=\{(1,2)\}$.
With $\kappa_{120}=1$, stabilizing controllers take the form:
\begin{equation}
\begin{aligned}
& u_1^\varepsilon(t,x) = a_1(x) + \tfrac{4\pi \sqrt{|a_{120}|}}\varepsilon \cos\left(\tfrac{2\pi t}{\varepsilon}\right),\\
& u_2^\varepsilon(t,x) = a_2(x) + \tfrac{4\pi \sqrt{|a_{120}|}}\varepsilon {\rm sign}\left(a_{120}\right) \cos\left(\tfrac{2\pi t}{\varepsilon}\right),
\end{aligned}
\label{Euler-feedback}
\end{equation}
$$
\begin{aligned}
\begin{pmatrix}
a_1(x) \\
a_2(x )\\
a_{120}(x)
\end{pmatrix} &= -{\mathcal F_2}^{-1}(x) \Bigl(\gamma x {+} f_0(x)\Bigr){=}{-}
\begin{pmatrix}
\gamma x_1{+}\alpha_1x_2x_3 \\
\gamma x_2{+}\alpha_2x_1x_3\\
\tfrac{\gamma}{2\alpha_3}x_1{+}\tfrac{1}{2}x_1x_2
\end{pmatrix}.
\end{aligned}
$$
Similarly to~\cite{ZuSIAM}, it can be shown that the trajectories of~\eqref{Euler} with controls~\eqref{Euler-feedback} converge to the origin \emph{exponentially}.

For the numerical simulation, we put
$
\alpha_1=3,\,\alpha_2=2,\,\alpha_3=1,\,\gamma=5,\,\varepsilon=1.
$
The resulting time-plots of $x_1(t)$, $x_2(t)$, $x_3(t)$  for two sets of initial conditions are presented in Fig.~\ref{fig_body} (in red). To illustrate other stabilizing strategies for a rotating rigid body with two control torques, we also show simulations results with smooth time-invariant controls from~\cite[p.~293]{Aeyels} (in green), discontinuous time-invariant controls from~\cite[Eq.~(22)]{Rey96} (in dark blue), and time-varying controls from~\cite[Eqs.~(24)--(25)]{Mor97} (in light blue). Note that the controls~\eqref{Euler-feedback} and~\cite[Eqs.~(24)--(25)]{Mor97} ensure exponential convergence to zero for the initial data from \emph{an entire neighborhood} of the origin, while the controls~\cite[Eq.~(22)]{Rey96} are not applicable if $x_1(0)=x_2(0)=0$, and the controls from~\cite[p.~293]{Aeyels} ensure only asymptotic (\emph{but not exponential}) convergence. The time-plots of all mentioned control functions are presented in Fig.~\ref{fig_cont}, left. The plot of the function $\tfrac{1}{t}\ln\|x(t)\|$  illustrates the exponential decay rate of the solutions of~\eqref{Euler} (see Fig.~\ref{fig_cont}, right).

\subsection{Stabilization of an underwater vehicle}
\begin{figure}[tpt]
\begin{center}
\includegraphics[width=1\linewidth]{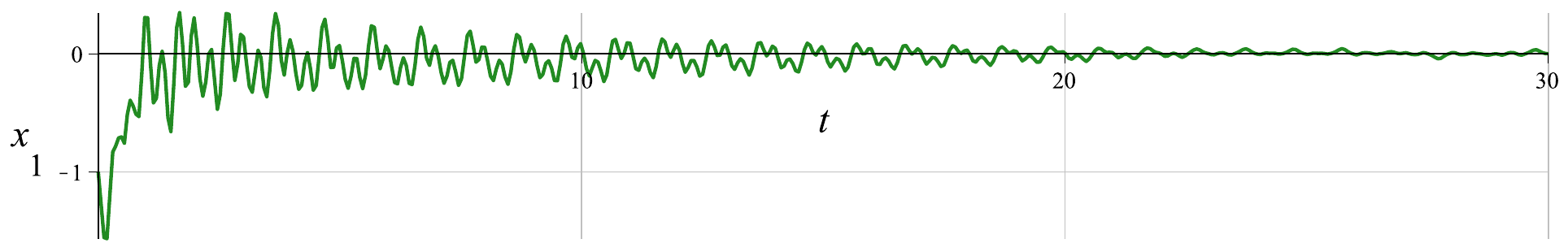}
\includegraphics[width=1\linewidth]{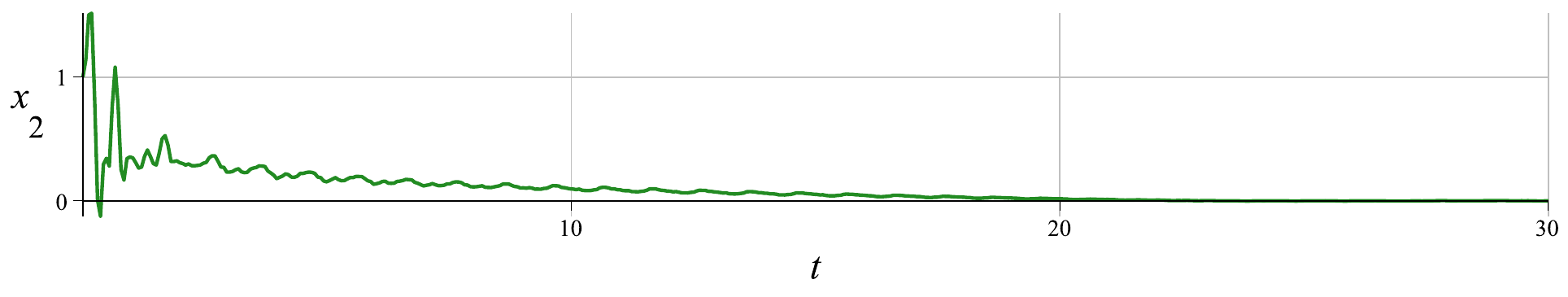}
\includegraphics[width=1\linewidth]{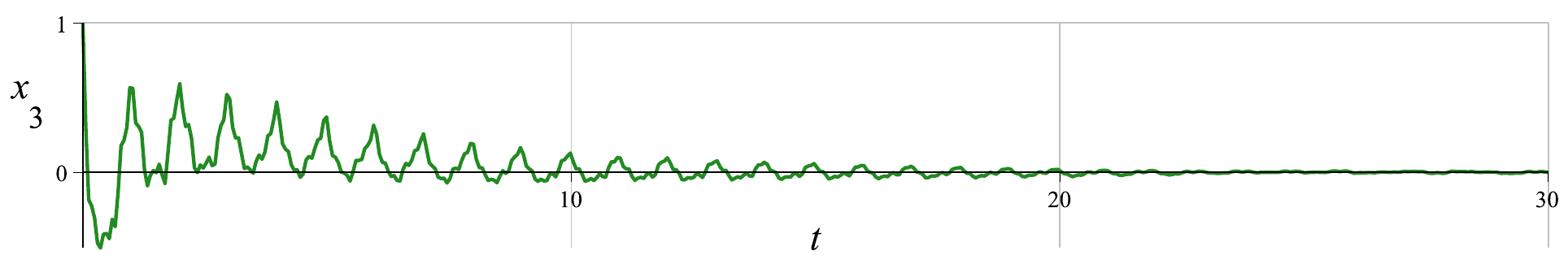}
\includegraphics[width=1\linewidth]{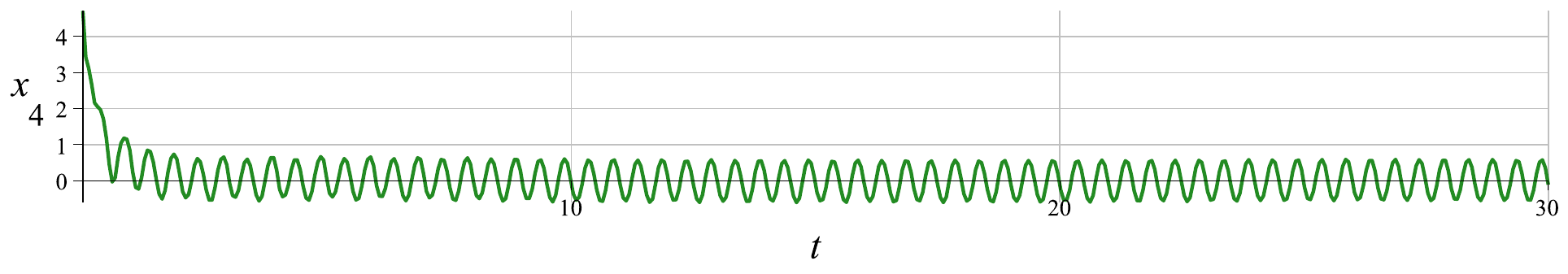}
\includegraphics[width=1\linewidth]{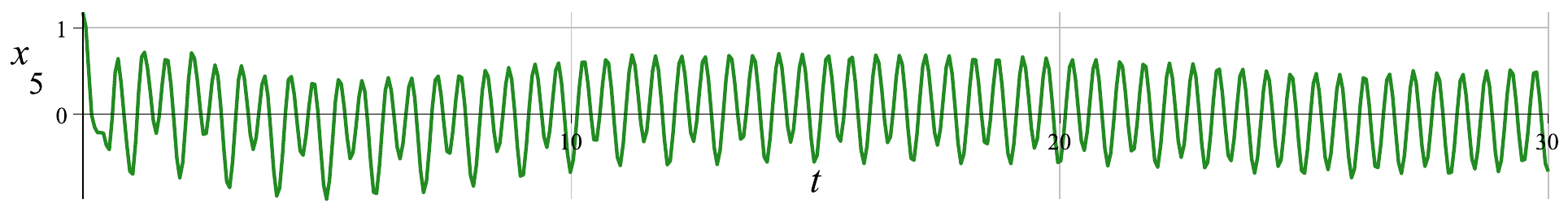}
\includegraphics[width=1\linewidth]{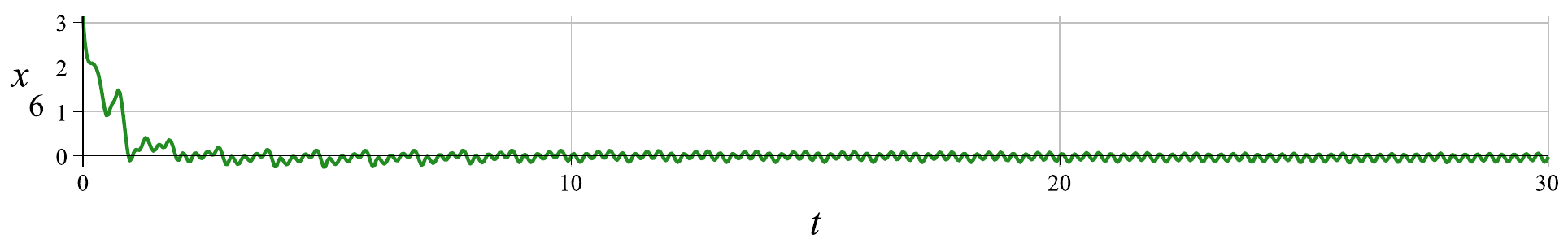}
\end{center}
\caption{Time-plots of the trajectories  for system~\eqref{underdr}--\eqref{cont_underdr} with the  initial condition
$x(0)=\big(-1,\ 1,\ 1,\ \tfrac{3\pi}{2},\ \tfrac{3\pi}{8},\ \pi\big)^\top$.
\label{fig_under}}
\end{figure}

In this subsection, we consider a nonlinear system with a non-vanishing drift term.
Namely, we consider the equations of motion of an autonomous  underwater vehicle studied, e.g., in~\cite{Bara}. As in~\cite{ZG17}, assume that one of the components of the angular velocity is uncontrolled and remains constant. Then the dynamics can be represented in the  control-affine form:
\begin{equation}\label{underdr}
\dot x=f_0(x)+\sum_{k=1}^3f_k(x)u_k,
 \end{equation}
where $(x_1, x_2, x_3)$  denote coordinates of the center of mass and $(x_4$, $x_5$, $x_6)$ specify the vehicle orientation,
$$
\begin{aligned}
&f_0(x)=\omega\big(0,0,0,\cos (x_4)\tan(x_5),{-}\sin (x_4),\cos (x_4)\sec (x_5)\big)^\top,\\
& f_1(x)=\big(\cos (x_5)\cos (x_6),\cos( x_5)\sin (x_6),-\sin (x_5),0,0,0\big)^\top,\\
& f_2(x) =\big(0,0,0,1,0,0\big)^\top, \\
& f_3(x)=\big(0, 0, 0, \sin (x_4)\tan(x_5), \cos(x_4),\sin (x_4)\sec (x_5)\big)^\top.
\end{aligned}
$$
The control $u_1$ represents the translational velocity  along the $x_1$-axis,  $(u_2,u_3)$ control the two components of the angular velocity, while the third component $\omega$ is assumed to be constant. Straightforward computations show that, for all $x\in D=\{x \in \mathbb R^6:-\tfrac{\pi}{2}<x_5<\tfrac{\pi}{2}\}$, system~\eqref{underdr} satisfies the  rank condition~\eqref{rank} with $S_1=\{1,2,3\}$, $S_{2}=\{(1,3),(2,3)\}$, $S_{10}=\{1\}$, $S_{20}=\emptyset$:
$
{\rm span}\Big\{f_1(x),\,f_2(x),\,f_3(x),\,[f_1,f_3](x),\,[f_2,f_3](x),\,[f_1,f_0](x)\Big\}=\mathbb R^n$ for all $x\in D.$
Following the control design algorithm described in Section~II.C, we take the following  functions:
\begin{align}
u_1^\varepsilon(t,x){=}a_1(x)&{+}2\sqrt{\tfrac{\pi \kappa_{13}|a_{13}(x)|}{\varepsilon}}{\rm sign}\big(a_{13}(x)\big)\cos\Big(\tfrac{2\pi\kappa_{13}t}{\varepsilon}\Big)\nonumber\\
&{+}\tfrac{2\pi\kappa_{10}a_{10}(x)}{\varepsilon}\sin\Big(\tfrac{2\pi\kappa_{01}t}{\varepsilon}\Big),\nonumber
\end{align}
\begin{align}
u_2^\varepsilon(t,x){=}a_2(x)&{+}2\sqrt{\tfrac{\pi \kappa_{23}|a_{23}(x)|}{\varepsilon}}{\rm sign}\big(a_{23}(x)\big)\cos\Big(\tfrac{2\pi\kappa_{23}t}{\varepsilon}\Big),\nonumber\\
u_3^\varepsilon(t,x){=}a_3(x)&{+}2\sqrt{\tfrac{\pi \kappa_{13}|a_{13}(x)|}{\varepsilon}}\sin\Big(\tfrac{2\pi\kappa_{13}t}{\varepsilon}\Big)\label{cont_underdr}\\
&{+}2\sqrt{\tfrac{\pi \kappa_{23}|a_{23}(x)|}{\varepsilon}}\sin\Big(\tfrac{2\pi\kappa_{23}t}{\varepsilon}\Big),\nonumber
\end{align}
where
$
a(x)=-\mathcal F^{-1}(x) \Bigl(\gamma x + f_0(x)\Bigr)
$
with $\gamma>0$,   $\mathcal F^{-1}(x)$ is inverse to
$
(
f_1(x)\ f_2(x)\ f_3(x)\ [f_1,f_3](x)\ [f_2,f_3](x)\ [f_1,f_0](x)
              ),
$
and $\kappa_{12},\kappa_{23},\kappa_{10}\in\mathbb N$ are pairwise distinct.
For the numerical simulation, we put
$
\omega=2,\,\kappa_{13}=1,\,\kappa_{23}=2,\,\kappa_{10}=3,\,\gamma=5,\,\varepsilon=1,
$
 Note that the drift term $f_0$ does not satisfy assumption~(A1), which means that $x^*=0\in\mathbb R^6$ is not an equilibrium of the corresponding closed-loop system.
However, as it is mentioned in Remark~\ref{rem_non0}, it is still possible to ensure practical convergence to $x^*$ (i.e. convergence to some $\Delta$-neighborhood of $x^*$), as it is demonstrated in Fig.~\ref{fig_under}.
In the future work, we expect to thoroughly analyze the relation between control parameters and the radius $\Delta$ of a neighborhood that can be achieved by the trajectories of the corresponding closed-loop system.
\section*{APPENDIX}
The explicit formulas for $\Omega_1,\Omega_2$ in~\eqref{volt_eps_1} and~\eqref{volt_eps_2} are given below:
$$
\begin{aligned}
\Omega_1(x^0,\varepsilon)&=\varepsilon\sum_{i\in S_1,l\in S_{10}}a_{i}(x^0)a_{l0}(x^0)\big[f_l,f_i\big](x^0)
\\
+&\tfrac{\varepsilon^2}{2}\Big\{ L_{f_0}f_0(x^0)+\sum_{i\in S_1}a_{i}(x^0)\Big(L_{f_0}f_{i}(x^0)+L_{f_{i}}f_0(x^0)\Big)
\\
+&\sum_{i_1,i_2\in S_1}L_{f_{i}}f_{i_2}(x^0)a_{i_1}(x^0)a_{i_2}(x^0)\Big\},
\end{aligned}
$$
 $$
\begin{aligned}
 \Omega_2&(x^0,\varepsilon)=\varepsilon\Big\{\sum_{(l_1,l_2)\in S_{20}}\big|a_{l_1l_20}(x^0)\big|\Big(\big[f_{l_1},[f_{l_1},f_0]\big](x^0)\\
 &+\big[f_{l_2},[f_{l_2},f_0]\big](x^0)\Big)+4\sum_{k_1,k_2=1}^m\sum_{\tiny \begin{array}{c}
                                                                        (l_1,l_2)\in S_{20} \\
                                                                        (j_1,j_2)\in S_{20} \\
                                                                        (l_1,l_2)\ne(j_1,j_2)
                                                                      \end{array}
 }\zeta_{k_1l_1l_2}\zeta_{k_2j_1j_2}\\
 &\times\tfrac{\kappa_{l_2l_20}\kappa_{j_1j_20}}{\kappa_{j_1j_20}^2-\kappa_{l_1l_20}^2}\sqrt{|a_{l_1l_20}||a_{j_1j_20}|}\Big(\big[f_{0},[f_{l_1},f_{l_2}]\big](x^0)\\
 &+\sum_{i\in S_1}a_{i}(x^0)\big[f_{i},[f_{l_1},f_{l_2}]\big](x^0)\Big)+\sum_{i\in S_1}\sum_{(l_1,l_2)\in S_{20}}a_i(x^0)\\
 &\times|a_{l_1l_20}(x^0)|\Big(\big[f_{l_1},[f_{l_1},f_{i}]\big](x^0)+\big[f_{l_2},[f_{l_2},f_{i}]\big](x^0)\\
{+}&{\rm sign(a_{l_1l_20}(x^0))}\big(\big[f_{l_1},[f_{l_2},f_{i}]\big](x^0){+}\big[f_{l_2},[f_{l_1},f_{i}]\big](x^0)\big)\Big)\Big\}
 \\
 &+\tfrac{\varepsilon^2}{2}\Big\{ L_{f_0}f_0(x^0)+\sum_{i\in S_1}a_{i}(x^0)\Big(L_{f_0}f_{i}(x^0)+L_{f_{i}}f_0(x^0)\Big)]]\\
&+\sum_{i_1,i_2\in S_1}L_{f_{i}}f_{i_2}(x^0)a_{i_1}(x^0)a_{i_2}(x^0)\Big\}\\
&+\tfrac{\varepsilon^2}{\pi}\Big\{
\sum_{k=1}^m\sum_{(l_1,l_2)\in S_{20}}\zeta_{kl_1l_2}\tfrac{\sqrt{|a_{l_1l_20}(x^0)|}}{\kappa_{l_1l_20}}\Big(\big[f_{0},[f_{0},f_{k}]\big](x^0)\\
&+\sum_{i\in S_1}a_i(x^0)\big(\big[f_{0},[f_{i},f_{k}]\big](x^0)+\big[f_{i},[f_{0},f_{k}]\big](x^0)\big)\\
&+\sum_{i_1,i_2\in S_1}a_{i_1}(x^0)a_{i_2}(x^0)\big[f_{i_1},[f_{i_2},f_{k}]\big](x^0)\Big)
\Big\}\\
&+\tfrac{\varepsilon^3}{6}\Big\{L_{f_0}L_{f_0}f_0(x^0)+\sum_{i\in S_1}a_i(x^0)\Big(L_{f_0}L_{f_0}f_i(x^0)\\
&+L_{f_0}L_{f_i}f_0(x^0)+L_{f_i}L_{f_0}f_0(x^0)\Big)+\sum_{i_1,i_2\in S_1}a_{i_1}(x^0)a_{i_2}(x^0)\\
&\times\Big(L_{f_0}L_{f_{i_1}}f_{i_2}(x^0)+L_{f_{i_1}}L_{f_0}f_{i_2}(x^0)+L_{f_{i_1}}L_{f_{i_2}}f_0(x^0)\Big)
\\
&+\sum_{i_1,i_2,i_3\in S_1}a_{i_1}(x^0)a_{i_2}(x^0)a_{i_3}(x^0)L_{f_{i_1}}L_{f_{i_2}}f_{i_3}(x^0)\Big\}.
\end{aligned}
 $$
 Here $\zeta_{kl_1l_2}=\delta_{kl_1}+\delta_{kl_2}{\rm sign}(a_{kl_1l_2}(x^0))$.
\end{document}